\documentclass[a4paper,12pt]{amsart}
\usepackage{amsfonts,amsthm,amssymb,amsmath,amscd}
\usepackage[colorlinks=true,linkcolor=blue,citecolor=blue]{hyperref}
\usepackage{amsmath}
\usepackage{tikz-cd}

\usepackage{tabularx}
\usepackage{tcolorbox}
\usepackage[english]{babel}
\usepackage{amssymb,amsmath}
\usepackage{xcolor}

\newtheorem{Theo}{Theorem}[section]
\newtheorem{Prop}[Theo]{Proposition}
\newtheorem{Coro}[Theo]{Corollary}

\newtheorem{Defi}[Theo]{Definition}

\newtheorem{Prob}[Theo]{Problem}

\newtheorem{Rema}[Theo]{Remark}

\newcommand{\Hcal}{\mathcal{H}}
\newcommand{\T}{\mathbb{T}}
\newcommand{\Bcal}{\mathcal{B}}
\newcommand{\D}{\mathbb{D}}
\newcommand{\C}{\mathbb{C}}
\newcommand{\Z}{\mathbb{Z}}

\def\N{\mathbb{ N}}
\def\R{\mathbb{ R}}

\begin{document}

\title{Hardy spaces of general Dirichlet series -- a survey}


\author[Defant]{Andreas Defant}
\address[]{Andreas Defant\newline  Institut f\"{u}r Mathematik,\newline Carl von Ossietzky Universit\"at,\newline
26111 Oldenburg, Germany.
}
\email{defant@mathematik.uni-oldenburg.de}

\author[Schoolmann]{Ingo Schoolmann}
\address[]{Ingo Schoolmann\newline  Institut f\"{u}r Mathematik,\newline Carl von Ossietzky Universit\"at,\newline
26111 Oldenburg, Germany.
}
\email{ingo.schoolmann@uni-oldenburg.de}

\maketitle

\begin{abstract}
\noindent
The  main purpose of this article is to  survey on some key elements of a recent  $\mathcal{H}_p$-theory of general
Dirichlet series $\sum a_n e^{-\lambda_{n}s}$, which was mainly inspired by the work of Bayart and Helson on ordinary Dirichlet series $\sum a_n n^{-s}$. In view of an ingenious identification of Bohr,
the $\mathcal{H}_p$-theory of ordinary Dirichlet series
  can be seen as a sub-theory of Fourier analysis on  the infinite dimensional
torus $\mathbb{T}^\infty$. Extending these ideas, the $\mathcal{H}_p$-theory of $\lambda$-Dirichlet series is
build as a sub-theory of Fourier analysis on what we call $\lambda$-Dirichlet groups. A number of problems is added.
\end{abstract}

\vspace{1cm}


\noindent
\renewcommand{\thefootnote}{\fnsymbol{footnote}}
\footnotetext{2010 \emph{Mathematics Subject Classification}: Primary} \footnotetext{\emph{Key words and phrases}:
} \footnotetext{}

\section{Introduction}

Within the last two decades the theory of  ordinary Dirichlet series $\sum a_{n} n^{-s}$ saw a sort of renaissance.
The study of these series in fact  was   one of the hot topics in mathematics at the beginning of the 20th. Among others,
H. Bohr, Besicovitch, Bohnenblust, Hardy,  Hille, Landau, Perron,  M.~Riesz, and Neder,
were  the leading mathematicians in this issue. The theory
 lived a sort of golden moment between the 1910’s and the 1930’s but after that it was somehow
forgotten. Some 20 years ago the seminal article \cite{HLS} of Hedenmalm, Lindqvist and Seip   called again the attention from analysis to Dirichlet
series. Since then a lot has been going on, and ordinary  Dirichlet series have been studied with new techniques from
functional and harmonic analysis.

 \medskip

Bohr's main interest was to derive properties of Dirichlet series from the analytical properties of
the holomorphic functions defined by them. It is well known that Dirichlet series converge on half planes,
and where they converge, they define a holomorphic function. Bohr considered three abscissas for a given
Dirichlet series $D = \sum a_n n^{-s}$: $\sigma_c(D)$, $ \sigma_u(D)$, and $\sigma_a(D)$ that define the maximal half planes on which the series respectively converges,
converges uniformly, or converges absolutely. He also considered a fourth abscissa, $\sigma_b(D)$, that gives the maximal
half plane on which the series defines a bounded and holomorphic function. Then
\[
\sigma_c(D) \leq \sigma_b(D) \leq \sigma_u(D) \leq \sigma_a(D)\,.
\]
In \cite{Bohr} Bohr was interested on describing the absolute convergence abscissa of an ordinary  Dirichlet series in terms of analytic properties of its limit function. He proved
\begin{equation} \label{sigma-ub}
\sigma_b(D) =  \sigma_u(D)\,,
\end{equation}
(see e.g. \cite[Corollary 1.14]{Defant}), and then considered the number
\[
S= \sup \{ \sigma_a(D) - \sigma_u(D) \colon D = \sum a_n n^{-s} \}\,,
\]
that gives the maximal width of the band on which a Dirichlet series can converge uniformly but not absolutely.
Bohr showed that $S \leq 1/2$. The problem of whether or not this was the correct value remained open
for some 15 years, until Bohnenblust and Hille in \cite{BoHi} indeed proved that
\begin{equation} \label{S12}
S=\dfrac{1}{2}
\end{equation}
(see also \cite[Theorem 4.1]{Defant}).
 These ideas are the seeds of a recent revival of interest in the research area opened up by these early contributions. A new field emerged intertwining the classical work in novel ways with modern functional analysis, infinite dimensional holomorphy, probability theory as well as analytic number theory. As a consequence, a number of  challenging research problems crystallized and were solved over the last decades. We refer to the monographs \cite{Helson},  \cite{QQ}, and \cite{Defant} were many   of the   key elements of this new developments
 for ordinary Dirichlet series are described in detail.

 \medskip

 A fundamental object in these investigations
 is given by the Banach space $\mathcal{D}_{\infty}$  of all ordinary Dirichlet series $D =\sum a_{n} n^{-s}$ which converge and define a bounded, and then necessarily holomorphic, function on $[Re>0]$
 (endowed with the supremum norm on $[Re>0]$), and a  celebrated result from   \cite{HLS}  shows that $\mathcal{D}_{\infty}$ in fact  equals the Hardy space $H_{\infty}(\T^{\infty})$ on the infinitely dimensional torus. Let  us explain this in more detail.
 The infinitely dimensional torus $\T^{\infty}$ is the  infinite product of $\mathbb{T} = \{w \in \mathbb{C}\colon |w|=1\}$ which forms a natural compact abelian group
 on which the Haar measure is given by the  normalized Lebesgue measure. The characters on $\mathbb{T}^\infty$, so the elements in the dual group, consists  of all  monomials $z \mapsto z^\alpha$, where
 $\alpha = (\alpha_k)\in \mathbb{Z}^{(\mathbb{N})}$ (all a finite sequence
 of integers), and $H_{\infty}(\T^{\infty})$ denotes the closed subspace of all $f \in L_\infty(\mathbb{T}^\infty)$ such that the Fourier coefficient
 \[
 \hat{f}(\alpha) = \int_{\T^{\infty}} f(w) w^{-\alpha} dw =0\,,
 \]
 whenever $\alpha < 0$ (in the sense that some $\alpha_k < 0$). Then, based on Bohr's work,  it is proven in \cite{HLS} (see again \cite[Corollary 5.3]{Defant}) that there is a unique linear isometry
 \begin{equation} \label{HLS}
  H_{\infty}(\T^{\infty}) = \mathcal{D}_{\infty}\,,
 \end{equation}
  which preserves Fourier- and Dirichlet coefficients in the sense that
 \begin{equation} \label{vision}
 \hat{f}(\alpha) = a_{n}\,\, \text{ whenever  } n= \mathfrak{p}^\alpha := p_1^{\alpha_1}\ldots  p_N^{\alpha_N}\,,
 \end{equation}
 where $\alpha = (\alpha_1, \ldots, \alpha_N, 0 \ldots) \in \mathbb{N}_0^{(\mathbb{N})}$ and $\mathfrak{p} = 2,3,5 \ldots$ the sequence of primes.
 We refer to this result as the  Bohr-Hedenmalm-Lindqvist-Seip theorem (see also \cite[Corollary 5.3]{Defant}).
 The crucial point of its proof  is the fact that each natural number has a unique prime number decomposition, as well as Kronecker's theorem on Diophantine approximation (see e.g.  \cite[Proposition 3.4]{Defant}):
 The continuous group homomorphism
\[
\beta: \mathbb{R}  \rightarrow \mathbb{T}^\infty\,, \,\,\,\, t \mapsto (p_k^{-it})_{k=1}^\infty
\]
has dense range,
which in particular implies that for each $n$ and $\alpha$ with  $n = \mathfrak{p}^\alpha$ the following diagram
commutes:
\begin{equation} \label{Kron}
\begin{tikzpicture}[scale = 0.8]
        \node (G) at (0,0) {$\mathbb{T}^\infty$};
        \node (T) at (3,0) {$\mathbb{T}$};
        \node (R) at (0,-2) {$\mathbb{R}$};

        \draw[-latex] (G) -- node[above] {$z^\alpha$} (T);
        \draw[-latex] (R) -- node[below right] {$e^{-it\log n}$} (T);
        \draw[-latex] (R) -- node[left] {$\beta$} (G);
    \end{tikzpicture}
    \end{equation}

    The original proof of \eqref{HLS} goes  a detour through infinite dimensional holomorphy using a result of Cole and Gamelin from \cite{ColeGamelin}. Denote by $H_\infty (B_{c_0})$ all holomorphic (Fr\'echet differentiable) functions
    $f$ on the open unit ball $B_{c_0}$ of $c_0$, which endowed with the sup norm forms a Banach space.  Then
    Cole and Gamelin show  that there  is a unique isometric linear bijection
    \begin{equation} \label{hol_p}
      H_\infty (\mathbb{T}^\infty) = H_\infty (B_{c_0})\,, ~f \mapsto g
    \end{equation}
which preserves Fourier and monomial coefficients in the sense that for every multiindex $\alpha$
\begin{equation}\label{FourMon}
\hat{f}(\alpha) = \frac{\partial^\alpha g(0)}{\alpha!}
\end{equation}
(see \cite[Theorem 5.1]{Defant}).
The equality of Banach spaces from \eqref{HLS} shows that a Dirichlet series $D$ belongs to $\mathcal{D}_\infty$ if and only if there is a function $f \in H_\infty(\mathbb{T}^\infty)$
    such that the Dirichlet coefficients $(a_n(D))$ and the Fourier coefficients $(\hat{f}(\alpha))$ coincide in the sense of \eqref{vision},
    if and only if there is a function $f \in H_\infty(B_{c_0})$
    such that  the Dirichlet coefficients $(a_n(D))$ and the monomial coefficients  $(\frac{\partial^\alpha g(0)}{\alpha!})$ coincide in the sense of \eqref{FourMon}. This links the theory of ordinary Dirichlet which generate bounded, holomorphic functions on the positive half plane intimately with  Fourier analysis on the group $\T^{\infty}$ as well as infinite dimensional holomorphy on the open unit ball of $c_0$.

\medskip

 More generally, Bayart in \cite{Bayart} developed an $H_p$-theory of Dirichlet series. Recall that the Hardy space $H_p(\mathbb{T}^\infty)\,, \, 1 \leq p \leq \infty,$ is the  closed subspace  of all $f \in L_p(\mathbb{T}^\infty)$
 such that $\hat{f}(\alpha) =0$ only if $\alpha < 0$. Then the  Banach spaces $\mathcal{H}_p$ of ordinary Dirichlet series by definition  is the isometric image of $H_p(\mathbb{T}^\infty)$ under the identification
 from \eqref{vision}:
 \[
 \mathcal{H}_p = H_p(\mathbb{T}^\infty)\,,
 \]
 and the Bohr-Hedenmalm-Lindqvist-Seip theorem  from  \eqref{hol_p}  reads
 \begin{equation} \label{fog}
 \mathcal{D}_\infty = \mathcal{H}_\infty = H_\infty(\mathbb{T}^\infty)\,.
 \end{equation}
 Again the Banach spaces $\mathcal{H}_p$  can be reformulated in terms of holomorphic functions in infinitely many variables.
 Instead of looking at holomorphic functions on $B_{c_{0}}$, we  look at holomorphic functions on
$\ell_{2} \cap B_{c_{0}}$, understood as an open subset of $ \ell_{2}$. Then
\[
H_{p}(\ell_{2} \cap B_{c_{0}})
\]
is the Banach space of all holomorphic functions $g: \ell_{2} \cap B_{c_{0}}\rightarrow \mathbb{C}$ for which
\begin{equation} \label{jungels}
\Vert g \Vert_{H_{p}(\ell_{2} \cap B_{c_{0}})} =
\sup_{n \in \mathbb{N}} \sup_{0<r<1} \left( \int_{\mathbb{T}^{n}} \vert g(rw_{1}, \ldots , r w_{n},0,0, \ldots) \vert^{p} d(w_{1}, \ldots , w_{n})
\right)^{\frac{1}{p}} < \infty  \,.
\end{equation}
Then for  $1\leq  p< \infty$ there is a unique isometric identity
\begin{equation} \label{Hphol}
H_{p}(\ell_{2} \cap B_{c_{0}}) = H_p(\mathbb{T}^\infty) \,.
\end{equation}
identifying Fourier coefficients and monomial coefficients in the sense of \eqref{FourMon}
(a detailed proof is given in \cite[Chapter 13]{Defant}).

\medskip

Let us  now face general Dirichlet series -- our main object of interest. Given  a frequency $\lambda:=(\lambda_{n})$ (i.e.  a non-negative sequence of real numbers tending to $+\infty$), a
$\lambda$-Dirichlet series is a formal sum
 $$\sum a_{n} e^{-\lambda_{n}s}$$
 with complex Dirichlet coefficients $(a_{n})$ and a complex variable $s$.
Denote by  $\mathcal{D}(\lambda)$  the space of all formal $\lambda$-Dirichlet series.

\medskip

We start recalling  some basic facts.
The natural domains of convergence of Dirichlet series are half spaces (see \cite[\S II.2, Theorem 1]{HardyRiesz}), and as in the ordinary case the 'abscissas'
\begin{align*}
&
\sigma_{c}(D)=\inf\left \{ \sigma \in \R \mid D \text{ converges on } [Re>\sigma] \right\},
\\&
\sigma_{a}(D)=\inf\left \{ \sigma\in \R \mid D \text{ converges absolutely on } [Re>\sigma] \right\},
\\&
\sigma_{u}(D)=\inf\left \{ \sigma \in \R \mid D \text{ converges uniformly on } [Re>\sigma] \right\},
\end{align*}
rule the convergence theory of general Dirichlet series.
  Again general Dirichlet series $D$ define holomorphic functions on $[Re>\sigma_{c}(D)]$, which relies on the fact that they converge uniformly on all compact subsets of $[Re>\sigma_{c}(D)]$ (see \cite[\S II.2, Theorem 2]{HardyRiesz}).

  \medskip

There are useful Bohr-Cahen formulas for the abscissas $\sigma_{c}$ and $\sigma_{a}$, that are, given $D=\sum a_{n}e^{-\lambda_{n}s}$,
 $$\sigma_{c}(D)\le \limsup_{N} \frac{\log\Big(\big| \sum_{n=1}^{N} a_{n}\big| \Big) }{\lambda_{N}} ~\text{ and  } ~\sigma_{a}(D)\le \limsup_{N} \frac{\log\Big(  \sum_{n=1}^{N} |a_{n}| \Big) }{\lambda_{N}},$$
 where in each case equality holds if the left hand side is non negative.
See \cite[\S II.6 and \S II.7]{HardyRiesz} for a proof. The formula for $\sigma_{u}$ (and its proof) extends from the ordinary case in \cite[\S 1.1, Proposition 1.6]{Defant} canonically to arbitrary $\lambda$'s:
\begin{equation} \label{sigmaU}
\sigma_{u}(D)\le \limsup_{N} \frac{\log\Big(\sup_{t \in \mathbb{R}} \big|\sum^{N}_{n=1} a_{n}e^{-\lambda_{n}it}\big|\Big)}{\lambda_{N}},
\end{equation}
where again equality holds if the left hand side is non negative.

 \medskip

We point out that making the jump from the ordinary case $\lambda=(\log n)$ to arbitrary frequencies reveals serious difficulties and that many of  the above  ideas fail for general Dirichlet series. Much of the ordinary theory relies on Bohr's theorem from \eqref{sigma-ub}, the fact that for each ordinary Dirichlet series the abscissa of uniform convergence and boundedness coincide. This phenomenon  fails  for general Dirichlet series. Further due to the prime number theorem each natural number $n$ has its prime number decomposition $n=\mathfrak{p}^{\alpha}$ and so  the frequency $(\log n)$ can be written as a linear combination of $(\log p_{j})$ with natural coefficients. This intimately links the theory of ordinary Dirichlet series  with the theory of holomorphic functions on polydiscs, and in particular with the theory of polynomials $\sum c_{\alpha} z^{\alpha}$
in finitely many complex variables. One of several consequences is that  $m$-homogeneous Dirichlet series $\sum a_n n^{-s}$, i.e.  $a_n \neq 0$ only if $n$ has $m$ prim factors, are linked with $m$-homogeneous polynomials.
 This way powerful tools enter the game, as e.g. polynomial inequalities (like the Bohnenblust-Hille inequalities, hypercontractivity of convolution with the Poisson kernel, etc.), $m$-linear forms, or  polarization. So modern Fourier analysis and infinite dimensional holomorphy enrich the theory of ordinary Dirichlet series considerably. All this is carefully explained in \cite{Defant}.

  \medskip

Unfortunately facing  general Dirichlet series many of these powerful bridges
 seem to collapse and new questions arise which make the theory of general Dirichlet quite  challenging.

  \medskip

 Nevertheless inspired by the ordinary theory  we are able to introduce a Fourier analysis setting for the study of  general Dirichlet series  by restricting ourserlves to $\lambda$-Dirichlet series with Dirichlet coefficients that  actually are Fourier coefficients defined by functions on compact abelian groups $G$ of a certain type; compact abelian groups allowing a  continuous homomorphism $\beta \colon \R \to G$ with dense range.

 \section{Bohr's  theorem}

To get started we define the space $\mathcal{D}_{\infty}(\lambda)$  analogously to the space $\mathcal{D}_{\infty} = \mathcal{D}_{\infty}((\log n))$.
\begin{Defi} \label{Dinfty} Let $\lambda$ be a frequency. Then $\mathcal{D}_{\infty}(\lambda)$ denotes  the space of all $\lambda$-Dirichlet series $D$ which converge on $[Re>0]$, and define
(a then necessarily holomorphic)   bounded
 function there.
\end{Defi}
 For instance, if $\lambda=(n)_{n=0}^\infty$, then looking at the transformation $z=e^{-s}$ we  easily conclude that  $\mathcal{D}_{\infty}((n))$ is simply $H_{\infty}(\D)$, the space of all bounded and holomorphic functions on the open unit ball $\D$.
  And if $\lambda=(\log n)$, then we are in the  ordinary case.

 \medskip

Recall that there is a unique coefficient preserving isometry identifying $H_{\infty}(\D)$ and $H_{\infty}(\T)$.
 Hence, in both cases, $\lambda = (n)$ and $\lambda = (\log n)$, the space $\mathcal{D}_{\infty}(\lambda)$ can be described in terms of Fourier analysis, that is, it can be considered as a Hardy space, namely $\mathcal{D}_{\infty}((n))=H_{\infty}(\T)$ and $\mathcal{D}_{\infty}((\log n))=H_{\infty}(\T^{\infty})$. In view of these two examples the following  question arises naturally:
 Given an arbitrary frequency $\lambda$,
is it possible to describe $\mathcal{D}_{\infty}(\lambda)$ in terms of a sort of Hardy space on a compact abelian group?

\medskip

The first step towards the solution of this problem, is to follow ideas of Bohr and Landau, who asked under which assumptions the abscissa $\sigma_u(D)$ can be described in terms of analytic properties of the limit function of $D$. Define
\[
\sigma_{b}(D)=\inf \sigma\,,
\]
where the infimum is taken over all $\sigma \in \R$ such that $D$ converges and defines a bounded function on $[\text{Re} > \sigma]$ (which then is automatically holomorphic).
Additionally define, provided $\sigma_{c}(D)<\infty$,
\begin{align*}
\sigma_{b}^{ext}(D)=\inf \sigma \,,
\end{align*}
the infimum now taken over all  $\sigma \in \R$ such that the limit function of $D$ allows a holomorphic and bounded  extension to $[\text{Re} > \sigma]$.
By definition $\sigma_{c}(D)\le \sigma_{b}(D) \le\sigma_{u}(D)\le \sigma_{a}(D)$ and $\sigma_{b}^{ext}(D)\le \sigma_{b}(D)$. But in  general all these abscissas differ. For instance an example of Bohr from \cite{Bohr4} shows that  $\sigma_{c}(D)=\sigma_{b}^{ext}(D)=\sigma_{b}(D)=-\infty$ and $\sigma_{u}(D)=+\infty$ is possible.

\medskip

A very prominent research project  at the beginning of  the 20th century was to find conditions for frequencies  $\lambda$
under which
\begin{equation} \label{problem}
\sigma_{b}^{ext}(D) = \sigma_{u}(D)
\end{equation}
holds for all somewhere convergent $\lambda$-Dirichlet series $D$.
In this context  the space $$\mathcal{D}_{\infty}^{ext}(\lambda)$$ of all somewhere convergent $\lambda$-Dirichlet series $D=\sum a_{n}e^{-\lambda_{n}s}$, which have a holomorphic and bounded extension to $[Re>0]$,
appears naturally. Endow
$\mathcal{D}_{\infty}^{ext}(\lambda)$  with the semi norm $\|D\|_{\infty}:=\sup_{[Re>0]}|f(s)|$, where $f$ is the (unique) extension of $D$. We will later in fact see that this always defines a norm --  but that,
in general, neither $\mathcal{D}_{\infty}(\lambda)$ nor $\mathcal{D}_{\infty}^{ext}(\lambda)$ form  Banach spaces.

\medskip

Much of the abstract theory of ordinary Dirichlet series is based on a fundamental theorem of Bohr from \cite{BohrStrip} which shows that  every $D \in \mathcal{D}^{ext}_\infty((\log n))$ converges uniformly on all half spaces $[\text{Re} > \varepsilon]$ for  all
$\varepsilon >0$, and this then easily implies that $\mathcal{D}_\infty^{\text{ext}}((\log n)) = \mathcal{D}_\infty((\log n))$ (see also \cite[Theorem 1.13]{Defant}).
For certain classes of frequencies $\lambda$,
   Bohr \cite{Bohr} and Landau \cite{Landau}  extended this result to $\lambda$-Dirichlet series.

\medskip

In \cite{Bohr} Bohr shows that $(\ref{problem})$ holds if $\lambda$ satisfies the following condition (we call it Bohr's condition $(BC)$):
\begin{equation} \label{BC}
\exists ~l = l (\lambda) >0 ~ \forall ~\delta >0 ~\exists ~C>0~\forall n \in \N: ~~\lambda_{n+1}-\lambda_{n}\ge Ce^{-(l+\delta)\lambda_{n}};
\end{equation}
roughly speaking this condition prevents the $\lambda_n$'s from getting too close too fast. Motivated by Bohr's work we introduce the following definition.

\begin{Defi}
We say that a frequency $\lambda$ satisfies Bohr's theorem, whenever every $D \in \mathcal{D}^{ext}_\infty(\lambda)$ converges uniformly on all half spaces $[\text{Re} > \varepsilon]\,,\,\,\varepsilon >0$, i.e. $\sigma_{b}^{ext}(D)= \sigma_{u}(D)$.
\end{Defi}

 Observe that $\lambda=(\log n)$ satisfies $(BC)$ with $l=1$ and so  Bohr's theorem holds.
In \cite{Landau} Landau gave a weaker sufficient condition than $(BC)$ (we call it Landau's condition $(LC)$), which extends the class of frequencies  which satisfy Bohr's theorem:
\begin{equation} \label{LC}
\forall ~\delta>0 ~\exists ~C>0~ \forall ~n \in \N \colon ~~ \lambda_{n+1}-\lambda_{n}\ge C e^{-e^{\delta\lambda_{n}}}.
\end{equation} We like to mention that in \cite[\S 1]{Neder} Neder went a step further and considered $\lambda$'s satisfying
$$\exists ~x>0 ~\exists ~C>0~ \forall ~n \in \N \colon ~~ \lambda_{n+1}-\lambda_{n}\ge C e^{-e^{x\lambda_{n}}}.$$
Then Neder proved that this condition is not sufficient for Bohr's theorem by constructing, given some $x>0$, a Dirichlet series $D$ (belonging to some frequency $\lambda$) for which $\sigma_{c}(D)=\sigma_{a}(D)=x$ and $\sigma_{b}^{ext}(D)\le 0$ holds, which in particular shows that $\mathcal{D}_{\infty}(\lambda) \subsetneq \mathcal{D}^{ext}_{\infty}(\lambda)$.

\medskip

Analysing Bohr's original proof from \cite{Bohr}, Bohr's theorem in the ordinary case was improved in \cite{BalasubramanianCaladoQueffelec} by a quantitative version (see again \cite[Theorem 1.13]{Defant}).
\begin{Theo} There is a constant $C>0$ such that for every $D \in \mathcal{D}_\infty((\log n))$ and  all  $N\ge 2$
\begin{equation} \label{BTordinary}
\sup_{[Re>0]} \bigg| \sum_{n=1}^{N} a_{n} n^{-s}\bigg| \le C \log(N) \|D\|_{\infty}.
\end{equation}
\end{Theo}
Note that in the ordinary case $\lambda=(\log n)$, this inequality together with $(\ref{sigmaU})$ implies $(\ref{problem})$.

\medskip

Like in the work of Bohr, Landau only proves the qualitative version of  the fact that each frequency $\lambda$ under his condition $(LC)$ satisfies
Bohr's theorem. To establish quantitative versions in the sense of  \eqref{BTordinary} means  to control the norm of the partial sum operator
$$S_{N} \colon \mathcal{D}_{\infty}^{ext}(\lambda) \to \mathcal{D}_{\infty}(\lambda), ~~ D \mapsto \sum_{n=1}^{N}a_{n}(D) e^{-\lambda_{n}s}\,.$$
 Using the summation method of typical means of order $k>0$ invented by M. Riesz (see Proposition~\ref{approx}), the following  estimate of $\|S_{N}\|$, which holds  under no further conditions on $\lambda$, is the main result of \cite[Theorem 3.2]{Schoolmann}.
\begin{Theo} \label{keylemma}
Let $D=\sum a_{n} e^{-\lambda_{n}s}\in \mathcal{D}_{\infty}^{ext}(\lambda)$. Then for all $0<k\le1$ and $N \in \N$:
$$\sup_{[Re>0]}\left| \sum_{n=1}^{N} a_{n}e^{-\lambda_{n}s} \right|\le C \frac{\Gamma(k+1)}{k} \left(\frac{\lambda_{N+1}}{\lambda_{N+1}-\lambda_{N}} \right)^{k} \|D\|_{\infty},$$
where $C>0$ is a universal constant and $\Gamma$ denotes the gamma function.
\end{Theo}
As a consequence, assuming Bohr's condition $(BC)$ for $\lambda$, the choice $k_{N}:=\frac{1}{\lambda_{N}}$, $N\ge 2$ (note that $\lambda_{1}=0$ is possible) leads to
\begin{equation} \label{coff}
\|S_{N}\colon \mathcal{D}_{\infty}^{ext}(\lambda) \to \mathcal{D}_{\infty}(\lambda)\|\le C_{1}(\lambda) \lambda_{N},
\end{equation}
which reproves $(\ref{BTordinary})$ for $\lambda=(\log n)$. Assuming Landau's condition $(LC)$,  the estimate from  Theorem~\ref{keylemma} with $k_{N}:=e^{-\delta \lambda_{N}}$, $\delta>0$, gives
\begin{equation} \label{ruecken}
\|S_{N}\colon \mathcal{D}_{\infty}^{ext}(\lambda) \to \mathcal{D}_{\infty}(\lambda)\|\le C_{2}(\lambda, \delta) e^{\delta\lambda_{N}};
\end{equation}
the quantitative version of Bohr's theorem under $(LC)$.

\medskip

 For particular frequencies $\lambda$ the bounds from \eqref{coff} and \eqref{ruecken}  may be bad, which isn't too surprising since these results in fact hold for {\it  all} frequencies  satisfying the conditions of Bohr or Landau. For instance,
 consider  the case $\lambda=(n)$. Then the projection
 \[
 S_N\colon H_{\infty}(\T)\to H_{\infty}(\T)\,,\,\,\,S_{N}(f)=\sum_{n=0}^{N} \hat{f}(n)z^{n}
 \]
  is nothing else than the convolution operator which assigns to every $f$ its convolution with the Dirichlet kernel $D_N$, and this immediately gives   that $\|S_{N}\|=\|D_{N}\|_{1}\sim \log(N)$. To our knowledge the optimal upper and lower bounds for the norm of $S_N$  in the ordinary case $\lambda=(\log(n))$ are still unknown.

  \begin{Prob} \label{P1} Determine optimal bounds for $\|S_{N}\colon \mathcal{D}_{\infty}((\log n)) \to \mathcal{D}_{\infty}((\log n))\|$.
\end{Prob}

  \medskip

The proof of Theorem \ref{keylemma} relies on the following independently interesting result
from \cite[Proposition 3.4]{Schoolmann} which was inspired by the work of Hardy and M.~Riesz
from \cite{HardyRiesz}.
\begin{Prop} \label{approx}
Let $D=\sum a_{n} e^{-\lambda_{n}s}\in D_{\infty}^{ext}(\lambda)$ with extension $f$. Then for all $k> 0$ the Dirichlet polynomials
$$R_{x}^{k}(D)=\sum_{\lambda_{n}<x} a_{n} \left(1-\frac{\lambda_{n}}{x}\right)^{k} e^{-\lambda_{n}s}$$
converge uniformly to $f$ on $[Re>\varepsilon]$ for all $\varepsilon>0$. Moreover,
\begin{equation} \label{stau}
\sup_{x\ge 0} \| R_{x}^{k}(D)\|_{\infty} \le  \frac{e}{2\pi}\frac{\Gamma(k+1)}{k}\|D\|_{\infty}.
\end{equation}
\end{Prop}
In the language of \cite{HardyRiesz} Proposition \ref{approx} states that, given any order $k>0$, then  on every  halfplane $[Re>\varepsilon]$ the limit function of a Dirichlet series $D\in \mathcal{D}^{ext}_{\infty}(\lambda)$ is  the uniform limit of its first typical means of order $k$.

\medskip

Moreover, Proposition \ref{approx} gives a direct link to the theory of almost periodic functions on $\R$, and proves that $\mathcal{D}_{\infty}^{ext}(\lambda)$ in fact is a normed space (see Corollary \ref{normed}). Note that a priori, $\|\cdot\|_{\infty}$ is only a semi norm, or equivalently, it is not clear whether
 $\mathcal{D}^{ext}_{\infty}(\lambda)$ can be considered as a subspace of $H_{\infty}[Re>0]$, the Banach space of all holomorphic and bounded functions on $[Re>0]$.
  Here it is important to distinguish Dirichlet series from their limit functions, and  to prove that $\|\cdot\|_{\infty}$  in fact is a norm on $\mathcal{D}^{ext}_{\infty}(\lambda)$  requires to check
that  all Dirichlet coefficients of $D$ vanish provided $\|f\|_\infty =0$.

\medskip

 Recall that by definition a continuous function $f \colon \R \to \C$ is called (uniformly) almost periodic,
 whenever for every $\varepsilon>0$ there is a number $l>0$ such that for all intervals $I\subset \R$ with $|I|=l$ there is a translation number $\tau \in I$ such that $\sup_{x\in \R} |f(x-\tau)-f(x)|\le \varepsilon$ (see \cite{Besicovitch} for more information). Then by a result of Bohr a bounded and continuous function $f$ is almost periodic if and only if it is the uniform limit of trigonometric polynomials on $\R$, which are of the form $p(t):=\sum_{n=1}^{N} a_{n}e^{-itx_{n}}$ for some  $x_1, \ldots,x_{N}\in \R$ (see e.g. \cite[\S 1.5.2.2, Theorem 1.5.5]{QQ}). In particular, the Dirichlet polynomials $R_{x}^{k}(D)$ stated in Proposition \ref{approx} considered as functions on vertical lines $[Re=\sigma]$ are almost periodic.
\begin{Coro} \label{almost periodic} If $D \in \mathcal{D}_{\infty}^{ext}(\lambda)$ with extension $f$, then the function $f_{\sigma}(t):=f(\sigma+it)\colon \R \to \C$ is almost periodic and
$$a_{n}(D)=\lim_{T\to \infty} \frac{1}{2T} \int_{-T}^{T} f(\sigma+it) e^{(\sigma+it)\lambda_{n}}$$ for all $\sigma >0$. In particular, $\sup_{n \in \N} |a_{n}|\le \|f\|_{\infty}$.
\end{Coro}

This result is taken from \cite[Corollary 3.8]{Schoolmann}, and the next corollary is then an immediate consequence.

 \begin{Coro} \label{normed}  $\mathcal{D}_{\infty}^{ext}(\lambda)$, and consequently also its subspace
 $\mathcal{D}_{\infty}(\lambda)$, are normed spaces for any frequency $\lambda$.
\end{Coro}

 Another particular consequence of  Corollary \ref{almost periodic} is as follows. Clearly, we  may deduce from this corollary that  $\|S_{N}\colon \mathcal{D}_{\infty}^{ext}(\lambda) \to \mathcal{D}_{\infty}(\lambda)\| \le N$
 for all $N$, so equality \eqref{problem} follows, i.e. $\lambda$ satisfies Bohr's theorem, whenever $$L(\lambda):=\limsup_{N \to \infty} \frac{\log(N)}{\lambda_{N}}=0.$$
We like to mention that the number $L(\lambda)$ has the following geometric meaning in terms of abscissas. In \cite[\S 3, Hilfssatz 3, Hilfssatz 2]{Bohr2} Bohr proved  that
\begin{equation} \label{Bohrforpresident}
L(\lambda)=\sigma_{c}\left(\sum e^{-\lambda_{n}s} \right)=
\sigma_{a}\left(\sum e^{-\lambda_{n}s} \right)
=
\sup_{D \in \mathcal{D}(\lambda)} \sigma_{a}(D)-\sigma_{c}(D).
\end{equation}
For instance  we have that  $L((n))=0$, and   see as a consequence that for power series we  up to $\varepsilon$ can't distinguished between uniform convergence and boundedness of the limit function.

\medskip

Besides $\lambda$'s with $(BC)$ or $(LC)$, there is another class of frequencies $\lambda$ for which Bohr's theorem holds. In \cite{Bohr5} Bohr proved that $\mathbb{Q}$-linearly independent frequencies $\lambda$ satisfy
$$\sigma_{b}^{ext}(D)=\sigma_{a}(D)$$
for all somewhere convergent $\lambda$-Dirichlet series $D$. The use of Kronecker's theorem, which states that the set $\left\{ (e^{-\lambda_{n}it})_n \mid t \in \R\right\}$ is dense in $\T^{\infty}$ whenever the real sequence $(\lambda_{n})$ is $\mathbb{Q}$-linearly independent,  combined with Proposition \ref{approx}
 in \cite[Theorem 4.7]{Schoolmann} lead to  an alternative proof of this fact.

\begin{Theo} \label{independent} Let $D=\sum a_{n} e^{-\lambda_{n}s} \in \mathcal{D}^{ext}_{\infty}(\lambda)$ with a $\mathbb{Q}$-linearly independent frequency $\lambda$. Then $(a_{n})\in \ell_{1}$ and $\|(a_{n})\|_{1}= \|D\|_{\infty}$. Moreover,
$$\mathcal{D}^{ext}_{\infty}(\lambda)=\mathcal{D}_{\infty}(\lambda)=\ell_{1}\,,\,\,\,\,\, \sum a_n e^{-\lambda_n s}  \mapsto (a_n)\,,$$
and $\lambda$ satisfies Bohr's theorem. In particular,
$$\sup_{N\in \N} \sup_{[Re>0]} \left| \sum_{n=1}^{N} a_{n} e^{-\lambda_{n}s} \right| =\|D\|_{\infty}.$$
\end{Theo}
The following theorem summarizes  some of the preceding results on frequencies $\lambda$ satisfying Bohr's theorem.

\begin{Theo} \label{Bohr theo}
A frequency $\lambda$ satisfies Bohr's theorem, if one of following  conditions  holds:
\begin{itemize}
\item
 $L(\lambda) =0$
  \item
   $\lambda$ is $\mathbb{Q}$-linearly independent
\item
$\lambda$  satisfies $(LC)$ (or the stronger condition $(BC)$)
 \end{itemize}
 Moreover, in each of these cases we have  $\mathcal{D}_\infty^{\text{ext}}(\lambda) = \mathcal{D}_\infty(\lambda)$, but non of these conditions is  necessary for Bohr's theorem.
\end{Theo}

Unfortunately $\mathcal{D}_{\infty}(\lambda)$ may fail to be a Banach space. The following result
from  \cite[Theorem 5.2]{Schoolmann} is inspired by a construction of Neder in \cite{Neder}.
\begin{Prop} \label{incomplete}
Let $\lambda$ be a frequency. Then there is a strictly increasing sequence $(s_{n})$ of natural numbers such that $\mathcal{D}_{\infty}(\eta)$, where $\eta$ is the frequency obtained by ordering the set
\begin{equation} \label{Neder}
\left\{ \lambda_{n}+\frac{j}{s_{n}}(\lambda_{n+1}-\lambda_{n}) \mid n \in \N, ~j=0, \ldots s_{n}-1\right\}
\end{equation}
increasingly, is not complete and $\mathcal{D}_{\infty}(\eta)\subsetneq \mathcal{D}_{\infty}^{ext}(\eta)$. In particular, Bohr's theorem fails for $\eta$.
\end{Prop}
The good news is that the completeness of $\mathcal{D}_{\infty}(\lambda)$ is guaranteed by  several sufficient conditions on $\lambda$ which in  concrete cases
have a good chance to be checked. The following collection of results was proved in \cite[Theorem 5.1]{Schoolmann} (for the case of $(BC)$ see also \cite{ChoiKimMaestre}).

\begin{Theo} \label{completeness}
$\mathcal{D}^{ext}_\infty(\lambda)$ is complete, if $L(\lambda)<\infty$. In particular, $\mathcal{D}_{\infty}(\lambda)$ is a Banach space and coincides with $\mathcal{D}_{\infty}^{ext}(\lambda)$ provided that
$\lambda$ satisfies one of the following conditions:
\begin{itemize}
\item
 $L(\lambda) =0$
  \item
   $\lambda$ is $\mathbb{Q}$-linearly independent
\item
$\lambda$  satisfies $(LC)$ and $L(\lambda)<\infty$ (this includes $(BC)$).
\end{itemize}
Moreover, none of these  conditions is necessary for the completeness of $\mathcal{D}_{\infty}(\lambda)$.
\end{Theo}
In view of the different nature of the stated conditions in Theorem \ref{Bohr theo} and Theorem \ref{completeness}, it seems that we are far away from a characterization of completeness of $\mathcal{D}_{\infty}(\lambda)$ (or $\mathcal{D}_{\infty}^{ext}(\lambda)$).
\begin{Prob} \label{P1} Characterize completeness of $\mathcal{D}_{\infty}(\lambda)$ and $\mathcal{D}_{\infty}^{ext}(\lambda)$.
\end{Prob}
In this context more  interesting questions appear naturally.

\begin{Prob} \label{P2}  Is there any relation between the  $\lambda$'s satisfying Bohr's theorem and the $\lambda$'s for which $\mathcal{D}_{\infty}(\lambda)$ is complete?
\end{Prob}
For instance the frequency $\lambda:=(\sqrt{\log n})$ fulfils $(LC)$ (and so satisfies Bohr's theorem), but we don't know wether  $\mathcal{D}_{\infty}((\sqrt{\log n}))$ is complete.
\begin{Prob}\label{P3}  Is $\mathcal{D}_{\infty}((\sqrt{\log n }))$ complete?
\end{Prob}
Recall Theorem \ref{Bohr theo} where we conditions on frequencies  under  which Bohr's theorem hold, that is $\sigma_{b}^{ext}(D)=\sigma_{b}(D)$
for all $\lambda$-Dirichlet series $D$.
\begin{Prob} \label{P4} Find a reasonable condition on the frequency $\lambda$  which
is  weaker than $(LC)$, but is sufficient for the equality  $\sigma_{b}^{ext}(D)=\sigma_{b}(D)$
for all $\lambda$-Dirichlet series $D$.
\end{Prob}
By Theorem \ref{completeness}, $\mathcal{D}_{\infty}(\lambda)$ is complete, if $\sigma_{b}^{ext}(D)=\sigma_{b}(D)$ for all $\lambda$-Dirichlet series $D$ and $L(\lambda)<\infty$. So any contribution to Problem \ref{P4} may give partial solutions to Problem \ref{P2}.

\medskip

We finish summarizing a few relations of the three conditions $(BC)$, $(LC)$ and  $L(\lambda)<\infty$
(see \cite[Remark 4.1]{Schoolmann}).

\begin{Rema} \label{connections} \text{ }
\begin{itemize}
\item $(BC)$ implies $L(\lambda)<\infty$ and $(LC)$.
\item $(LC)$ plus $L(\lambda)<\infty$ does not necessarily imply  $(BC)$.
\item $L(\lambda)<\infty$ does not necessarily imply  $(LC)$, and so neither $(BC)$.
\item $(LC)$ does not necessarily imply   $L(\lambda)<\infty$, and so neither $(BC)$.
\end{itemize}
\end{Rema}

\section{Hardy spaces of general Dirichlet series} \label{Hardysection}
Now we start our $\mathcal{H}_{p}$-theory on general Dirichlet series. As already mentioned in the introduction we restrict ourselves to general Dirichlet series with Dirichlet coefficients which actually are Fourier coefficients of  functions on  certain compact
abelian groups.
This has several  advantages.
One is  that  the class of all general Dirichlet series simply is  too large to obtain a good understanding.
 Assuming that the Dirchlet coefficients are Fourier coefficients gives more structure and allows to use tools from harmonic analysis like  the Hausdorff-Young inequality or Plancherel's theorem (among others).
   A further advantage of our setting is that  Bayart's  $\mathcal{H}_{p}$-theory of ordinary Dirichlet series embeds in a natural way. Whereas the $\mathcal{H}_{p}$-theory of ordinary Dirichlet series is basically Fourier
  analysis on the infinite dimensional torus  $\T^{\infty}$, this group  fails to be the right model
  for general Dirichlet series. In fact, the Bohr compactification $\overline{\R}$ of $\R$ and products of $\widehat{\mathbb{Q}_d}$ (the dual group of the rationals endowed with the discrete topology) turn out to be  suitable substitutes. Finally, fixing some $\lambda$, regarding the different realisations of $\lambda$-Dirichlet series of this type, another feature of our approach is that the $\mathcal{H}_{p}$-theory of general Dirichlet series we intend to present will be independent of the chosen suitable group for $\lambda$.

\medskip

We are interested in the subclass of all  compact abelian groups $G$ which allow a continuous homomorphism $\beta\colon \R  \to G$ with dense range. We call such  pairs $(G,\beta)$  Dirichlet groups.
Hence, given such a pair,
  the characters $x = e^{-ix\cdot}  \in  \widehat{\beta}(\widehat{G}) \subset \widehat{\mathbb{R}}$ are precisely those for which there is a unique character  $h_{x}  \in \widehat{G}$ with  $e^{-ix\cdot} = h_x \circ \beta$;
to understand this recall that $\mathbb{R} = \hat{\mathbb{R} }, x \to e^{-ix\cdot} $ is a group isomorphism.
      In particular, we have
\begin{equation} \label{dualmapping3}
\widehat{G}=\{h_{x} \mid x \in \widehat{\beta}(\widehat{G}) \}.
\end{equation}
Then the following notion (first given in \cite{DefantSchoolmann}) turns out to be fundamental for our purposes.
\begin{Defi}
Let $\lambda$ be a frequency and $(G, \beta)$  a Dirichlet group. Then $G$ is called a $\lambda$-Dirichlet group whenever $\lambda \subset \widehat{\beta}(\widehat{G})$, which
 means that  the following diagram
   commutes for every $n\in \N$:
 \begin{equation*}
\begin{tikzpicture}[scale = 0.8]
        \node (G) at (0,0) {$G$};
        \node (T) at (3,0) {$\mathbb{T}$};
        \node (R) at (0,-2) {$\mathbb{R}$};

        \draw[-latex] (G) -- node[above] {$h_{\lambda_n}$} (T);
        \draw[-latex] (R) -- node[below right] {$e^{-i\lambda_n \cdot}$} (T);
        \draw[-latex] (R) -- node[left] {$\beta$} (G);
    \end{tikzpicture}
    \end{equation*}
 \end{Defi}

Let us give examples. Denoting by $\mathfrak{p}:=(p_{n})$ the sequence of prime numbers, the compact group $\T^{\infty}$ together with the mapping
$$\beta_{\T^{\infty}} \colon \R \to \T^{\infty}, ~~ t \mapsto \mathfrak{p}^{-it}$$
forms a $(\log n)$-Dirichlet group (see again Kronecker's theorem). This example keeps us in track to recover results on ordinary Dirichlet series. The 'mother' of all possible examples is as follows:
Given a subgroup $U$ of $\R$, the topological group  $\widehat{(U,d)}$ together with the
 mapping
$$ \beta_{\widehat{(U,d)}}\colon \R \to  \widehat{(U,d)}, ~~ t \mapsto \left[u \mapsto e^{-itu}\right]$$
 forms a  Dirichlet group. So in particular, for $U = \mathbb{Z}$ and identifying $\T =\widehat{\mathbb{Z}}$, we
 obtain the $(n)$-Dirichlet group $(\T, \beta_{\T})$, where
\[
\beta_{\T}: \R \rightarrow \T, \, t \mapsto e^{-it}.
\]
The compact abelian group  $\overline{\R}:=\widehat{(\R,d)}$ is the so-called   Bohr compactification
of $\R$ which forms a $\lambda$-Dirichlet group for all frequencies $\lambda$ with the embedding
\begin{equation*}
\beta_{\overline{\R}}\colon (\R,  |\cdot|) \hookrightarrow \overline{\R},\,\, x \mapsto \left[ t \mapsto e^{-ixt}\right]\,.
\end{equation*}
Besides $\overline{\R}$, countable products of $\widehat{\mathbb{Q}_d}:=\widehat{\left(\mathbb{Q}, d\right)}$, where $d$ denotes the discrete topology, also form 'universal' Dirichlet groups. To explain this , let $B:=(b_{1}, b_{2}, \ldots)$ be a $\mathbb{Q}$-linearly independent sequence of real numbers of length $N \in \mathbb{N} \cup \{\infty\}$.
Then
\begin{equation} \label{themap}
T_B: \bigoplus_{n=1}^{N} \mathbb{Q} \hookrightarrow \R, ~~\alpha \mapsto \sum \alpha_{j}b_{j}
\end{equation}
is an  injective homomorphism, and hence its dual map
$$\widehat{T_{B}} \colon \R \to \widehat{\bigoplus_{n=1}^{\infty} \mathbb{Q}_{d}}, ~~ t \mapsto \left[ (q_{j})_{j} \mapsto e^{-it\sum q_{j}b_{j}}\right]$$
 has dense range. Since $\prod_{n=1}^{N} \widehat{\mathbb{Q}_d}=\widehat{\bigoplus_{n=1}^{\infty} \mathbb{Q}_{d}}$, the pair $\left(\prod_{n=1}^{N} \widehat{\mathbb{Q}_d}, \widehat{T_{B}}\right)$ is a Dirichlet group, which is 'universal' in the following sense: For every frequency $\lambda$ there is a suitable $B$ such that $\left(\prod_{n=1}^{N} \widehat{\mathbb{Q}_d}, \widehat{T_{B}}\right)$ is a $\lambda$-Dirichet group.

 \medskip

 To see this, let us recall that to every $\lambda$ there is a $\mathbb{Q}$-linearly independent sequence $B:=(b_{n})$ of real numbers, called basis for $\lambda$ (which can always be chosen as a subsequence of $\lambda$), such that $\lambda_{n}=\sum r^{n}_{k} b_{n}$ for some (unique) finite rational sequence $(r^{n}_{k})$. In this case, $R:=(r^{n}_{k})_{n,k}$ is said to be a Bohr matrix of $\lambda$ with respect to the basis $B$ and we write  $\lambda=(R,B)$. Hence $\left(\prod_{n=1}^{N} \widehat{\mathbb{Q}_d}, \widehat{T_{B}}\right)$ is a $\lambda$-Dirichlet group for every $\lambda$ with decomposition $\lambda=(R,B)$.

 \medskip

 Given a  $\lambda$-Dirichlet group $(G,\beta)$ and $1\le p \le \infty$, we define the Banach space
$$
H_{p}^{\lambda}\left(G\right):= \left\{ f \in L_{p}(G) \mid
\text{ $\hat{f}: \widehat{G} \to \mathbb{C}$ is supported by all $h_{\lambda_n}, n \in \mathbb{N}$}\right\}\,,$$
 and use it to define the following natural scale of Hardy spaces of general Dirichlet series.

\begin{Defi} \label{zugcrash}
The Hardy space $\Hcal_{p}(\lambda)$ of $\lambda$-Dirichlet series
is the  space
of all  $\sum a_n e^{-\lambda_n s}$ for which there is some
$f \in H_p^\lambda(G)$ such that   $a_n = \widehat{f}(h_{\lambda_{n}})$ for all $n$.
\end{Defi}

Together with
the norm $\|D\|_{p}:=\|f\|_{p}$ the space $\Hcal_{p}(\lambda)$ clearly forms a Banach space, and then by definition the Bohr map
\begin{equation*}
\Bcal \colon H^{\lambda}_{p}(G) \hookrightarrow \mathcal{D}(\lambda),
\,\,\, f \sim \sum_{\gamma \in \widehat{G}} \hat{f}(\gamma) \gamma\,\,\, \mapsto\,\,\, \sum_{n \in \N} \widehat{f}(h_{\lambda_{n}}) e^{-\lambda_{n}s}
\end{equation*}
 gives an isometric onto isomorphism
\[
\Hcal_{p}(\lambda) = H_{p}^{\lambda}(G)\,.
\]
The following fact, proved in \cite[Theorem 3.19]{DefantSchoolmann}, is fundamental.

\begin{Theo} \label{independentofgroup}
 $\Hcal_{p}(\lambda)$ is independent of the chosen $\lambda$-Dirichlet group $G$.
\end{Theo}
So the above definition of $\Hcal_{p}(\lambda)$  actually coincides with Bayart's definition of $\Hcal_{p}$
from \cite{Bayart} in the ordinary case (see also \cite[\S 11]{Defant}). Moreover, recall that in this case the groups $\T^{\infty}$ and $\overline{\R}$ are suitable $(\log n)$-Dirichlet groups, which immediately leads to the following consequence extending \eqref{HLS}.
\begin{Coro} \label{oorrddii} For all $1\le p \le \infty$ we have
$$H_{p}^{\log(n)}(\overline{\R})=\mathcal{H}_{p}(\log(n))=H_{p}(\T^{\infty}),~~ \|f\|_{p}=\|D\|_{p}=\|g\|_{p},$$
where $\widehat{f}(h_{\log n})=a_{n}(D)=\widehat{g}(\alpha)$, if $n=\mathfrak{p}^{\alpha}.$
\end{Coro}

Let us mention that there also is an internal description of $\Hcal_{p}(\lambda)$ through $\lambda$-Dirichlet polynomials without considering $\lambda$-Dirichlet groups. We denote by $Pol(\lambda)$ the space of all $\lambda$-Dirichlet polynomials $D(s)=\sum_{n=1}^{N} a_{n} e^{-\lambda_{n}s}$. For such polynomials we define
\begin{equation*} \label{limitHp}
\|D\|_{p}:=\left(\lim_{T \to \infty} \frac{1}{2T} \int_{-T}^{T}\left|\sum^{N}_{n=1} a_{n}e^{-\lambda_{n}it} \right|^{p} dt\right)^{\frac{1}{p}}.
\end{equation*}
Then this limit exists and gives a norm on $Pol(\lambda)$ (see e.g. \cite[Theorem 1.5.6]{QQ}
or \cite[Theorem 11.9]{Defant}). The following description was noted in \cite[Theorem 3.25]{DefantSchoolmann}.

\begin{Theo} \label{comp_pol}Let $1\le p < \infty$ and $\lambda$ be a frequency. Then the space $\Hcal_{p}(\lambda)$ is the completion of $(Pol(\lambda), \|\cdot\|_{p})$.
\end{Theo}
\medskip
\subsection{Frequencies of integer type}
Recall that given a frequency $\lambda$ there is a basis $B$ for $\lambda$ and a Bohr matrix $R$ such that $\lambda=(R,B)$. Let us look at the  ordinary case $\lambda=(\log n)$. Then $\log n= \sum \alpha_{j} \log p_{j}$
whenever $n=\mathfrak{p}^{\alpha}$, and hence a basis is given by $B=(\log p_{j})$ and every multi index $\alpha >0$ appears as a row in the corresponding Bohr matrix $R$. Moreover recall that $\T^{\infty}$ is a $(\log n)$-Dirichlet group with $\beta: \mathbb{R} \to \T^{\infty},\, \beta(t):=\mathfrak{p}^{-it}$.

\medskip

More generally, we call of frequency $\lambda$ of integer (natural) type, if there is a basis $B=(b_{n})$ such that the Bohr matrix $R$ associated to $B$ and $\lambda$ only has integer (natural) entries. In this case $\T^{\infty}$ is a $\lambda$-Dirichlet group with $\beta(t):=e^{-itB}:=(e^{-itb_{1}}, e^{-itb_{2}}, \ldots)$, $t \in \R$, and the rows of the corresponding  Bohr matrix are multi indices $\alpha$ (with integer entries). Actually one can show that a frequency $\lambda$ is of integer type if and only if there is a homomorphism $\beta \colon \R \to \T^{\infty}$ such that $(\T^{\infty},\beta)$ is a $\lambda$-Dirichlet group (see \cite[Remark 3.32]{DefantSchoolmann}).

\medskip

To see an example, consider the set $M:=\{n+\sqrt{2}m \mid n,m \in \mathbb{Z}\}$ which is dense in $\R$. Then any frequency $\lambda \subset M$ is of integer type with basis $B=(1,\sqrt{2})$. See \cite[Example 3.36]{DefantSchoolmann} for  a frequency which is not of integer type.

\medskip

Let $B$ be a basis of length $N \in \N\cup \{\infty\}$ for some frequency $\lambda$ of integer type with Bohr matrix $R$, and let us  write $\alpha \in R$ if the multi index $\alpha$ appears as a row in $R$. Then we define $H_{p}^{R}(\T^N)$ to be the space of all $g \in L_{p}(\T^{N})$ for which $\widehat{g}(\alpha)\ne 0$ implies $\alpha \in R$. By Theorem \ref{independentofgroup}, the $\mathcal{H}_{p}(\lambda)$'s can be identified with $H_{p}^{R}(\T^{N})$ in the following sense (see Theorem \cite[Theorem 3.30]{DefantSchoolmann}).
\begin{Theo} \label{integertype} Let $1\le p \le \infty$ and  $\lambda=(R,B)$  a frequency of integer type with basis of length $N \in \N\cup\{\infty\}$. Then there is a unique onto isometry $$\psi \colon \Hcal_{p}(\lambda) \to H_{p}^{R}(\T^N),  ~~ D \mapsto g$$ such that $\widehat{g}( \alpha)=a_{n}(D)$ for all multi indices $\alpha \in  R$ and $\lambda_{n}=\sum \alpha_{j} b_{j}$.
\end{Theo}
An immediate consequence is that the ordinary $\mathcal{H}_{p}$'s in the following sense are the largest spaces for $\lambda$'s of natural type.
\begin{Coro} \label{ordinary} Let $1\le p \le \infty$ and  $\lambda=(R,B)$  a frequency of natural type. Then there is a unique into isometry
$$ \psi \colon \Hcal_{p}(\lambda) \hookrightarrow \Hcal_{p}, ~~ \sum a_{n} e^{-\lambda_{n}s} \to \sum b_{n} n^{-s}$$
such that $a_{n}=b_{\mathfrak{p}^{\alpha}}$ for all $\alpha \in R$, where $\lambda_{n}=\sum \alpha_{j} b_{j}$.
\end{Coro}
\section{Some structure theory}

We summarize further properties of $\mathcal{H}_{p}(\lambda)$ which extend important key stones from Bayart's theory of ordinary Dirichlet series to our new theory
of general Dirichlet series.

\medskip
\subsection{Coincidence of $\pmb{H_\infty}$'s}
Note first that there are now two `$H_\infty$-spaces of $\lambda$-Dirichlet series' around, namely $\mathcal{D}_{\infty}(\lambda)$ and $\Hcal_{\infty}(\lambda)$. Recall  that by Proposition \ref{incomplete} there are  frequencies $\lambda$
such that $\mathcal{D}_{\infty}(\lambda)$ is not complete, hence in these cases  $\mathcal{D}_{\infty}(\lambda)\neq\mathcal{H}_{\infty}(\lambda)$. The following result is
 given in   \cite[Theorem 4.10, 4.12]{DefantSchoolmann}, and an far reaching extension of the Bohr-Hedenmalm-Lindqvist-Seip theorem from \eqref{fog}.

\begin{Theo} \label{Dstar} For every frequency $\lambda$ there is a coefficient preserving injective contraction
\begin{equation*} \label{Psi}
\mathcal{D}^{ext}_{\infty}(\lambda) \subset \mathcal{H}_\infty(\lambda)\,,
\end{equation*}
and if $\mathcal{D}^{ext}_{\infty}(\lambda)= \mathcal{D}_{\infty}(\lambda)$ and  $L(\lambda)<\infty$
(see again Theorem~\ref{Bohr theo}), then the embedding is even an isometric equality
\begin{equation*} \label{Psi2}
\mathcal{D}_{\infty}(\lambda) = \mathcal{H}_\infty(\lambda)\,.
\end{equation*}
\end{Theo}

We mention an interesting by-product of Theorem~\ref{Dstar} which is an immediate consequence of the  definition of $\mathcal{H}_2(\lambda)$ and  Parseval's equality.

\begin{Coro}
\label{H2}
For each  $D \in \mathcal{D}^{ext}_{\infty}(\lambda)$ we have  $(a_{n}(D))\in  \ell_{2}$ with $\|(a_{n}(D))\|_{2}\le \|D\|_{\infty}$, i.e. the embedding
\begin{equation*} \label{Psi3}
\mathcal{D}^{ext}_{\infty}(\lambda) \subset \mathcal{H}_2(\lambda)\,,
\end{equation*}
is a well-defined contraction.
\end{Coro}

\medskip

\subsection{Schauder bases} \label{Schauderbasischapter}
Given a frequency $\lambda$ and $1 < p < \infty$,
the following question is fundamental:
Do the $e^{-\lambda_{n} s}$ form a Schauder basis of $\Hcal_{p}(\lambda)$? In the ordinary case $\lambda =(\log n)$
the answer is affirmative as discovered in \cite{Saksman}, and we will see that the same result holds true whenever $\lambda$ is an arbitrary frequency.

\medskip

In \eqref{dualmapping3} we note that, given a Dirichlet group $(G, \beta)$, we have
$$\widehat{G}=\{h_{x} \mid x \in \widehat{\beta}(\widehat{G}) \}\,.$$
 This in particular shows that the dual group $\widehat{G}$ inherits the order of $\mathbb{R}$,
and hence we  deduce from \cite[Theorem 8.7.2.]{Rudin62}, that the 'Riesz projection'
 \begin{equation*}
\Phi\left(\sum a_{k} h_{x_{k}}\right):=\sum_{x_{k}\ge 0} a_{k} h_{x_{k}}
\end{equation*}
 is bounded on the subspace $Pol(G)$ of all polynomials in $\mathcal{H}_p(G)$, $1<p<\infty$. Then  standard arguments shows the following important theorem
 from \cite[Theorem 4.16]{DefantSchoolmann}.

\begin{Theo} \label{schauderbasisinHp}
Let $1 < p < \infty$ and $\lambda$ be a frequency. Then the monomials $e^{-\lambda_{n} s}$ form a Schauder basis for $\Hcal_{p}(\lambda)$.
\end{Theo}

An equivalent formulation of  Theorem \ref{schauderbasisinHp}  is that
for $1<p<\infty$
all projections
\[
S_N^p: \mathcal{H}_p(\lambda) \to \mathcal{H}_p(\lambda),\,\,
\sum a_n e^{-\lambda_n s} \mapsto \sum_{n=1}^N a_n e^{-\lambda_n s}
\]
 are uniformly bounded. But for the border cases  $p=1$ and $p=\infty$ this in general is false
 (e.g., for the frequencies  $\lambda = (\log n)$  or $\lambda= (n)$).
 Upper bounds for the growth of the partial sum operators in $\Hcal_{\infty}(\lambda)$
 were given in \eqref{coff} and \eqref{ruecken}. The following result handles the case $p=1$, and its proof in a sense  reduces to the case $p=\infty$ (see \cite[Proposition 4.17]{DefantSchoolmann}
 and \cite[12.5]{Defant} in the ordinary case).

 \begin{Prop} \label{basisconstant} Let $\lambda$ be a frequency. Assuming $(BC)$ for $\lambda$ there is a constant $C=C(\lambda)$ such that for all $N\ge 2$
\[
\|S_N^1: \mathcal{H}_1(\lambda) \to \mathcal{H}_1(\lambda)\|\leq C \lambda_{N}\,,
\]
and, assuming $(LC)$ and $L(\lambda)<\infty$, for every $\delta$ there is a constant $D=D(\delta,\lambda)$ such that  for all $N$
\[
\|S_N^1: \mathcal{H}_1(\lambda) \to \mathcal{H}_1(\lambda)\| \le De^{\delta\lambda_{N}}\,.
\]
\end{Prop}
\medskip
\subsection{Brothers Riesz theorem}
Recall the classical brothers Riesz theorem with states that the Hardy space $\mathcal{H}_{1}((n))=H_{1}(\T)$ coincides with the space $M(\T)$ of all bounded, regular and analytic Borel measures on $\T$. The corresponding result in ordinary case reads $\mathcal{H}_{1}((\log n))=H_{1}(\T^{\infty})=M(\T^{\infty})$ (due to Helson and  Lowdenslager from \cite{HelsonLowdenslager}). See also \cite[Theorem 13.5]{Defant} for a proof within the setting of ordinary Dirichlet series.  The brothers Riesz theorem extends to the case of general Dirichlet series.

\begin{Theo} \label{brotherRiesz}
Let $\lambda$ be any frequency and let $G$  a $\lambda$-Dirichlet group. Then the map
$$H_{1}^{\lambda}(G) \to M_{\lambda}(G), ~~f \mapsto f ~dm$$
 is an onto isometry. In particular $\mathcal{H}_{1}(\lambda)=M_{\lambda}(G)$.
\end{Theo}
Actually it  was discovered in \cite[Theorem 4.25]{DefantSchoolmann} that this is a fairly simple consequence of a more general result due to
\cite[Theorem 4]{Doss}; see \cite[\S 4.7]{DefantSchoolmann} for a discussion on this.

\medskip
\subsection{Montel  theorem}
In \cite[Lemma 18]{Bayart} Bayart proved that for every bounded sequence $(D^{N}) \subset \mathcal{D}_{\infty}((\log n))$ there is  a subsequence $(D^{N_{k}})$ and  $D \in \mathcal{D}_{\infty}((\log n))$ such that $(D^{N_{k}})$ converges to $D$ on $[Re> \varepsilon]$ for all $\varepsilon>0$. This fact,  sometimes called 'Montel's theorem', extends to the following classes of $\lambda$'s. See also \cite{QQ}, and \cite[Theorem 3.11]{Defant},
where Bayart's Montel theorem is deduced from a  Montel  type theorem for $H_\infty(B_{c_0})$
combined with \eqref{hol_p} and \eqref{fog}.

\begin{Theo}\label{montel} Let $\lambda$ satisfy $L(\lambda)=0$,  or $L(\lambda)<\infty$ and $(LC)$, or let $\lambda$ be $\mathbb{Q}$-linearly independent. Then for every $1\le p \le \infty$ and every
 bounded sequence $(D^{N}) \subset \Hcal_{p}(\lambda)$ there is  a  subsequence $(D^{N_{k}})$ and $D \in \Hcal_{p}(\lambda)$ such that  for all $\varepsilon>0$ the translations about $\varepsilon$
$$D^{N_{k}}_{\varepsilon}=\sum a_{n}(D^{N_{k}})e^{-\varepsilon \lambda_{n}} e^{-\lambda_{n}s}$$ converge to $D$ in $\Hcal_{p}(\lambda)$ as $k$ tends to $\infty$.
\end{Theo}

The proof in a sense reduces the $\Hcal_{p}(\lambda)$-case to the $\Hcal_{\infty}(\lambda)$-case.
This reduction   modifies an idea from \cite{AntonioDefant} showing  that, under the
assumption of the completeness of   $\mathcal{D}_{\infty}(\lambda)$, the map
$$ \Phi \colon \Hcal_{p}(\lambda) \hookrightarrow \mathcal{D}_{\infty}(\lambda, \Hcal_{p}(\lambda)), ~~ \sum a_{n} e^{-\lambda_{n}s} \mapsto \sum \left(a_{n} e^{-\lambda_{n}s}\right)  e^{-\lambda_{n}z}$$
is an into isometry; here $\mathcal{D}_{\infty}(\lambda, \Hcal_{p}(\lambda))$ stands for
the Banach space of all $\lambda$-Dirichlet series $\sum A_n e^{-\lambda_n s}$ with coefficients
$A_n \in \mathcal{H}_p(\lambda)$, which converge and are bounded on $[\text{Re}>0]$.

\medskip

\subsection{Hilbert's criterion}
As already mentioned the 'ordinary' $\mathcal{D}_{\infty}$ isometrically  equals $H_{\infty}(B_{c_{0}})$,
identifying Dirichlet and monomial coefficients. A crucial argument in the proof of this result (see \cite[\S 2.3]{Defant}) is that a continuous function $f \colon B_{c_{0}} \to \C$ is in $H_{\infty}(B_{c_{0}})$ if and only if all restriction maps $f_{N} \colon \D^{N}\to \C$ are in $H_{\infty}(\D^{N})$ and $\sup_{N} \|f_{N}\|_{\infty}<\infty$.

\medskip

Formulated for  ordinary Dirichlet series this shows that  $D=\sum a_n n^{-s}\in \mathcal{D}_{\infty}$ if and only if
all of its  so-called $N$-th abschnitte
$
D|_{N} = \sum a_{n} n^{-s}\,,
$
where the sum is only taken over those $n$ which have only the first $N$ primes as divisors,
belong to  $\mathcal{D}_{\infty}$ with uniformly bounded norms.
   In more vague terms,  $D \in \mathcal{D}_\infty$ if and only if all its finite dimensional blocks are in
   $\mathcal{D}_\infty$ with uniformly bounded norms (see also \cite[Theorem 3.11]{Defant}).

   \medskip

This phenomenon is also true in the general case. Given some decomposition $\lambda=(R,B)$ the $N$-th abschnitt $D|_{N}$ of a $\lambda$-Dirichlet series $D$ is the sum $\sum a_{n}(D) e^{-\lambda_{n}s}$, where $a_{n}(D)\ne 0$ implies  that $\lambda_{n}$ depends only on the first $N$ basis elements $b_{1}, \ldots, b_{N}$.

\medskip

 As proven in \cite[Theorem 4.22]{DefantSchoolmann}, general  Dirichlet series $D\in \mathcal{H}_{p}(\lambda)$ under
appropriate assumptions on the frequency   are again 'determined by their finite dimensional parts'.
\begin{Theo} \label{Nabschnitt}
Let $1\le p \le \infty$ and let $\lambda=(R,B)$ be a frequency satisfying one of the conditions of Theorem \ref{montel}. Let $D$ be a formal $\lambda$-Dirichlet series. Then $D \in \Hcal_{p}(\lambda)$ if and only if the $N$-th abschnitt $D|_{N} \in \Hcal_{p}(\lambda)$ for all $N \in \N$ and $\sup_{N} \|D|_{N}\|_{p}<\infty$. Moreover, in this case $\|D\|_{p}=\sup_{N \in \N} \|D|_{N}\|_{p}.$
\end{Theo}
By Theorem \ref{Dstar} we obtain the following particular case: If $\lambda$ satisfies $\sigma_{b}^{ext}=\sigma_{b}$ and $L(\lambda)<\infty$, then $D\in \mathcal{D}_{\infty}(\lambda)$ if and only if $D|_{N} \in \mathcal{D}_{\infty}(\lambda)$ for all $N \in \N$ and $\sup_{N} \|D|_{N}\|_{\infty}<\infty$.
\medskip
\subsection{Helson's theorem} \label{Helson}
A celebrated result of
Helson \cite{Helson3}
 on general Dirichlet series $\sum a_{n} e^{-\lambda_{n}s}$ states that if $\lambda$  satisfies $(BC)$ and $(a_{n})$ is $2$-summable, then for almost all homomorphism $\omega \colon (\R,+) \to \T$ the Dirichlet series $\sum a_{n} \omega(\lambda_n)e^{-\lambda_{n}s}$
converges on the open right half plane $[Re>0]$, or equivalently, for all $u>0$ the series
$\sum a_{n} e^{-\lambda_{n}u} h_{\lambda_n}$ converges almost everywhere on $\overline{\mathbb{R}}$
(note that if we for all rational $u$'s collect zero sets in
$\overline{\mathbb{R}}$, then we get the  aforementioned result).
 See also
\cite[Theorem 4.4]{HLS} for the ordinary case.

\medskip

Using our teminology, let  us rephrase this result for square summable functions on Dirichlet groups: Given a frequency $\lambda$ with $(BC)$, a $\lambda$-Dirichlet group $(G,\beta)$, and some  $f \in H_2^\lambda(G)$, then for every $u >0$ the series
\[
\sum_n \hat{f}(h_{\lambda_n}) e^{-\lambda_n u} h_{\lambda_n}
\]
converges almost everywhere in $G$.

\medskip

This result extends to functions
in $H_p^\lambda(G)$. To see this recall that for $u >0$ the Poisson kernel
$$P_{u}(t):=\frac{1}{\pi}\frac{u}{u^{2}+t^{2}} \colon \mathbb{R} \to \mathbb{R}$$
is integrable with $\|P_{u}\|_{1}=1$. Hence, by the Riesz representation theorem, for every $\lambda$-Dirichlet group $(G,\beta)$ the  functional
$$T_{u} \colon C(G)\to \mathbb{C}, ~~ g \mapsto \int_{\R} (g \circ \beta)(t) P_{u}(t) dt$$
defines a  measure  on $G$ which we call the 'Poisson measure' on $G$  denoted by $p_{u}$.

\medskip

The following theorem from \cite{DefantSchoolmann2} is a far reaching extension of Helson's theorem. It in fact shows that, given an appropriate frequency
$\lambda$, Helson's theorem extends to functions $f$ from our
scale $H_p^\lambda(G), 1 \leq p < \infty$, modelled on $\lambda$-Dirichlet groups $G$. Moreover, it describes the $G$-a.e. pointwise limit of
the Fourier series  $\sum_{n=1}^{\infty} \widehat{f}(h_{\lambda_{n}})e^{-\lambda_{n}u} h_{\lambda_{n}}$
in terms of convolution of $f$ with the Poisson measure $p_u$ on $G$, and -- most important -- it adds the relevant maximal inequality.

\begin{Theo}\label{maximalineq}
Let $(G,\beta)$ be a $\lambda$-Dirichlet group for  a frequency  with $(BC)$, and  $f \in H_p^\lambda(G)$
with $1 \leq p < \infty$. Then for every $u >0$
\[
f \ast p_{u} = \sum_{n=1}^{\infty} \widehat{f}(h_{\lambda_{n}})e^{-\lambda_{n}u} h_{\lambda_{n}}
\]
almost everywhere on $G$.
 Moreover, there is a constant $C=C(u, \lambda)$ such that for all $f \in H_{p}^{\lambda}(G)$
 \begin{equation} \label{ineqmaxi}
  \bigg\|\sup_{N} \Big| \sum_{n=1}^{N}\widehat{f}(h_{\lambda_{n}})e^{-\lambda_{n}u} h_{\lambda_{n}} \Big|
  \bigg\|_p \le C\|f\|_{p}.
 \end{equation}
\end{Theo}

\medskip

 \subsection{Vertical Limits}
  Let us  reformulate Theorem \ref{maximalineq} in terms of general Dirichlet series. Given a  $\lambda$-Dirichlet series $D(s)=\sum a_{n} e^{-\lambda_{n}s}$ and a $\lambda$-Dirichlet group $G$, we for $\omega \in G$ call the Dirichlet series
$$D^{\omega}(s)=\sum a_{n} h_{\lambda_n}(\omega) e^{-\lambda_{n}s}$$
a vertical limit of $D$. The following characterization of vertical limits justifies their name
(\cite[Proposition 4.6]{DefantSchoolmann}).

\begin{Rema} Let $(G, \beta)$ be a $\lambda$-Dirichlet group, and $D=\sum a_{n}e^{-\lambda_{n}s}$ a $\lambda$-Dirichlet series with $\sigma_{a}(D)\le 0$. Then
\begin{itemize}
\item[(1)]
For every  $\omega \in G$ there is a sequence $(\tau_{k}) \subset \R$ such  for all $\varepsilon>0$ the 'vertical translations'
 $$D_{i \tau_k} =  \sum a_n e^{-i \lambda_n \tau_k}e^{-\lambda_n s}$$
 converge to $D^{\omega}$ uniformly on $[Re>\varepsilon]$.

\item[(2)]
 Assume conversely that the $(D_{i\tau_k})$ for  some $(\tau_{k}) \subset \R$ converge
 uniformly on $[Re>\varepsilon]$
 to a holomorphic function $f$ for every  $\varepsilon>0$. Then there is $\omega \in G$ such that
    $f(s)= D^{\omega}(s)$
    for all $s \in [Re>0]$\,.
\end{itemize}
\end{Rema}

Note that the notion of vertical limits allows to rephrase Helson's theorem from the preceding section as follows: Given the $\lambda$-Dirichlet group $G=\overline{\R}$ and   $D\in \mathcal{H}_{2}(\lambda)$, then   $\overline{\R}$-almost all vertical limits $D^{\omega}$ converge on $[Re>0]$.

\medskip

The following result is then a straight forward reformulation of Theorem~\ref{maximalineq}.

\begin{Theo}\label{maximalineq2} Let $(G,\beta)$ be a $\lambda$-Dirichlet group for  a frequency  with $(BC)$, and  $D \in \mathcal{H}_p(\lambda)$
with $1 \leq p < \infty$. Then
 almost all vertical limits
 $$D^{\omega}(s) = \sum a_{n}(D) h_{\lambda_n}(\omega) e^{-\lambda_{n}s},\, \omega \in G,$$
 converge on $[Re >0]$. Moreover,
 for all $u>0$ there is a constant $C=C(u, \lambda)$ such that for all $D \in \Hcal_{p}(\lambda)$
 \begin{equation} \label{ineqmaxiDirichlet}
 \bigg(\int_{G} \sup_{N} \Big| \sum_{n=1}^{N} a_{n} h_{\lambda_n}(\omega) e^{-u\lambda_{n}} \Big|^{p} dm(\omega)\bigg)^{\frac{1}{p}} \le C\|D\|_{p}.
 \end{equation}
\end{Theo}

\begin{Rema} \label{ptoinfity}
Note that in the maximal inequality $(\ref{ineqmaxiDirichlet})$ the  constant $C = C(u, \lambda$)  is independent of $p$. Hence, if $p \to \infty$, then we recover that frequencies with $(BC)$ satisfy Bohr's theorem, i.e. for every somewhere convergent $\lambda$-Dirichlet series $D$ we have  $\sigma_{b}^{ext}(D)=\sigma_{u}(D)$.
\end{Rema}

Again convolution allows to relate the values of the vertical limits $D^{\omega}$ of a Dirichlet series $D\in \mathcal{H}_{p}(\lambda)$ and its associated function $f\in H_{p}^{\lambda}(G)$.
The vertical limits $D^{\omega}$ in some sense recover the associated function $f$ on vertical lines.
Note that, given  $f\in L_{p}(G)$,  almost all $\omega \in G$ define a measurable function $f_{\omega}(t):=f(\omega \beta(t))$, which itself is defined for almost all $t\in \R$. Hence, the convolution $$f_{\omega}*P_{u}(t):=\int_{\R} f_{\omega}(t-y) P_{u}(y) dy$$ is defined almost everywhere on $\R$.

\begin{Coro}\label{helsontheo}
Let  $1\le p < \infty$, $\lambda$ satisfy $(BC)$, and $(G,\beta)$ be a $\lambda$-Dirichlet group. If $f \in H_p^\lambda(G)$ and  $D\in \mathcal{H}_p(\lambda)$ are associated to each other,
then for all $u>0$ and almost all $t \in \R$
\begin{equation*} \label{formula}
D^{\omega}(u+it)=(f_{\omega}*P_{u})(t).
\end{equation*}
\end{Coro}

Let us  consider the particular case of ordinary Dirichlet series. We denote by  $\Xi$ the set of all completely multiplicative characters $\chi\colon \N \to \T$
(that is $\chi(nm)=\chi(n)\chi(m)$ for all $m$,$n$)
which with the pointwise multiplication forms an abelian  group. Looking at the group isomorphism
\begin{align} \label{chi1}
\iota\colon \Xi \to \T^{\infty}, ~~\chi \mapsto (\chi(p_{n}))_{n},
\end{align}
 where $\mathfrak{p}=(p_{n})$ again denotes the sequence of prime numbers, we see that $\Xi$ also  forms a compact abelian group, and its Haar measure $d \chi$ is the push forward measure of $dz$  through $\iota^{-1}$. Then applying Theorem \ref{maximalineq} to the  $(\log n)$-Dirichlet group $\Xi$ we obtain the following consequence.
\begin{Coro} \label{BayartsHelsonsthm}
Let $1\le p <\infty$. Then for all $u>0$ there is a constant $C=C(u)$ such that for all
Dirichlet series $\sum a_n n^{-s} \in \mathcal{H}_p(\lambda)$
$$\left(\int_{\Xi} \sup_{N} \left| \sum_{n=1}^N a_n \chi(n)n^{-u}\right|^{p} d\chi \right)^{\frac{1}{p}} \le C\|D\|_{p}\,.$$
In particular, almost all vertical limits $D^{\chi}(s) = \sum  a_n \chi(n)n^{-s}$, $ \chi \in \Xi$,  converge on $[Re >0]$.
\end{Coro}

\begin{Rema}
Helson in \cite{Helson} deduces from this result the curious fact that Riemann's conjecture holds for almost all
randomized $\zeta$-functions $\sum_n \chi(n) n^{-s}$, in the sense that almost all of them  converge
pointwise on $[\text{Re} >1/2]$ without having any zero in this half plane (see also \cite{HLS}).
\end{Rema}
Bayart in \cite{Bayart} deduces the 'almost everywhere part' of Corollary \ref{BayartsHelsonsthm} for the case $p=2$ from the Menchoff-Rademacher theorem on almost everywhere convergence of orthonormal series, and then he extends it to the range $1\le p<\infty$ using the 'hypercontractivity of the Poisson kernel'.

\medskip

Let us explain what is meant  by this. Recall the definition of the Poisson kernel on $\mathbb{T}$
$$p_{r}(w):=\sum_{n=-\infty}^{\infty} r^{|n|} w^{n}, ~~ 0<r<1,$$
and consider for all $1\le p<q <\infty$ the following convolution operator
$$\mathcal{P}_{r}\colon H_{p}(\T) \to H_{q}(\T), ~~f \mapsto f*p_{r}(\omega)=\sum_{n=0}^{\infty}\widehat{f}(n)(rw)^{n}.$$
An easy computation shows that $\|P_{r}\|\le \frac{1}{1-r}$ (independent of $p,q$), but this bound improves considerably depending on the relation of $p,q,r$. Weissler in \cite{Weissler} shows that $\|\mathcal{P}_{r}\|=1$ if and only if $r^{2}\le \frac{p}{q}$; a result which is refered to as the 'hypercontractivity of the Poisson kernel'.
Then using this result 'in each variable separately' shows that for all $u>0$ and $1\le p,q< \infty$ the translation operator
$$\mathcal{P}_{u}\colon \mathcal{H}_{p} \mapsto \mathcal{H}_{q}, ~~ D \mapsto D_{u}=\sum a_{n}(D) n^{-u} n^{-s}$$ is
well-defined and bounded (see \cite[Theorem 12.9]{Defant} for the details on all that). This explains, why in the ordinary case 'the almost everywhere part' of Corollary \ref{BayartsHelsonsthm} follows in the range $1\le p <\infty$,  once this result is proven for $p=2$.

\medskip

It seems that in the case of general Dirichlet series no 'hypercontractivity results' are available, and so our proof of Theorem \ref{maximalineq2} in \cite{DefantSchoolmann2} follows an entirely  different approach.
Nevertheless the following problem seems interesting in itself.

\begin{Prob} Given   $D = \sum a_n  e^{\lambda_n s}\in \mathcal{H}_{p}(\lambda), 1 \leq p < \infty$. Under which assumptions on $\lambda$
is it true  that  $D_{u} = \sum a_n e^{\lambda_n u} e^{\lambda_n s}\in \mathcal{H}_{q}(\lambda)$ for all $1 \leq q < \infty$ and   $u>0$?
\end{Prob}

From \cite[\S 12]{Defant} we know that for $1 \leq p < \infty$
 \begin{equation} \label{sec12}
\sup_{D \in \mathcal{H}_p} \sigma_c(D)= \sup_{D \in \mathcal{H}_p} \sigma_a(D) = \frac{1}{2}
\end{equation}
 The second equality was implicitly proved in \cite{Bayart}, and in view of the coefficient preserving identity  $\mathcal{H}_p =  H_p(\mathbb{T}^\infty)$, it implies
that the infimum over all $u > 0$ such that for each $f \in H_p(\mathbb{T}^\infty)$
the series $$\sum_{\alpha} \hat{f}(\alpha)\,\Big(\frac{z}{\mathfrak{p}^{u}}\Big)^\alpha$$ converges (here we of course talk about absolutele convergence) equals
$ 1/2$ (compare also with \eqref{Hphol}). The next result shows that there is no such 'wall' if we order the sum 'blockwise according to the order of $\mathbb{N}$'.

\begin{Coro}
Let  $f \in H_p(\mathbb{T}^\infty), 1 \leq  p < \infty$.  Then  for all $u >0$
\[
\lim_N \sum_{\mathfrak{p}^\alpha \leq N} \hat{f}(\alpha)\,\Big(\frac{z}{\mathfrak{p}^{u}}\Big)^\alpha \,\,\,\,\, \text{a.e. on $\mathbb{T}^\infty$}\,,
\]
and moreover
\[
\bigg( \int_{\mathbb{T}^\infty}   \sup_N \Big| \sum_{\mathfrak{p}^\alpha \leq N} c_\alpha
\Big(\frac{w}{\mathfrak{p}^{u}}\Big)^\alpha  \Big|^p  d \omega \bigg)^{1/p}
 \leq C(u) \|f\|_p\,,
\]
where $C= C(u)$ only depends on $u$.
\end{Coro}
What about the case $\omega=1$ (the neutral element of $G$) in Theorem \ref{maximalineq2}, i.e. (pointwise) convergence of Dirichlet series from $\mathcal{H}_{p}(\lambda)$? By $(\ref{Bohrforpresident})$ we always have $\sigma_{c}(D)\le   \sigma_{a}(D)\leq L(\lambda)$ for all $D \in \mathcal{H}_{p}(\lambda)$
(since in this case the sequence $(a_n(D))$ of Dirichlet coefficients is always bounded). In the particular case $\lambda=(\log n)$,
we know from \eqref{sec12} that
$$\sup_{D \in \mathcal{H}_{p}} \sigma_{c}(D)= \sup_{D \in \mathcal{H}_{p}} \sigma_{a}(D)=\frac{L((\log n))}{2}.$$
\begin{Prob} Given a frequency $\lambda$, give reasonable upper and lower estimates for
$\sup_{D \in \mathcal{H}_{p}(\lambda)} \sigma_{c}(D)$
and
$\sup_{D \in \mathcal{H}_{p}(\lambda)} \sigma_{c}(D)$\,,
for example in terms of $L(\lambda)$.
\end{Prob}

\medskip
\subsection{Carleson's theorem} \label{Carlesonchapter} What happens, if we consider the case $u=0$ in the preceding two sections, in particular in Theorem~\ref{maximalineq}?

\medskip

It turns out that this question is intimately related with one of the most
famous theorems in Fourier analysis, namely Carleson's theorem: The Fourier series of every $f \in L_2(\mathbb{T})$
converges almost everywhere. More generally, the Carleson-Hunt theorem states that if $f \in L_{p}(\T)$, $1<p<\infty$, then $f(z)=\sum_n \widehat{f}(n)z^{n}$ for almost all $z \in \T$ and
\[
  \Big\| \sup_N \big| \sum_{n=-N}^N \hat{f}(n)  z^n  \big| \Big\|_p
  \leq C \|f\|_p\,,
 \]
Recall that a frequency is of integer type if it has a  basis $B$ for which the Bohr matrix $R$ only contains integers.
For $\lambda$'s of this type the following theorem extends the  Carleson-Hunt theorem to our setting of $\lambda$-Dirichlet groups, and for such $\lambda$'s it is a strong improvement
of Theorem~\ref{schauderbasisinHp}.

 \begin{Theo} \label{Carleson}
Let $\lambda$ be a frequency of integer type,  $(G, \beta)$ a $\lambda$-Dirichlet group, and
  $1 < p < \infty$. Then for every $f \in H_p^\lambda(G)$
  $$f = \sum_{n=1}^\infty  \hat{f}(h_{\lambda_n})  h_{\lambda_n}  $$
  almost everywhere on $G$, and
   \[
  \bigg\| \sup_N \Big| \sum_{n=1}^N \hat{f}(h_{\lambda_n})  h_{\lambda_n}  \Big| \bigg\|_p
  \leq C \|f\|_p\,,
 \]
where $C = C(p)$ is a constant which only depends on $p$.
  \end{Theo}

  The case $p=2$ is due to  \cite[Theorem 1.4]{HedenmalmSaksman}, whereas the proof of the more general case
  given here follows closely the ideas of  this article combining them with the  Carleson-Hunt theorem and a magic trick of Fefferman from \cite{Fefferman}. This  is carried out in \cite{DefantSchoolmann2}.

  \medskip

It is somewhat surprising that Theorem \ref{Carleson} doesn't require $(BC)$ for $\lambda$.
\begin{Prob} Does Theorem \ref{maximalineq} hold for a wider class of frequencies than the class of $\lambda$'s satisfying $(BC)$? Is it even true without any condition on $\lambda$?
\end{Prob}

    Theorem \ref{Carleson}  easily  translates into a sort of Carleson-Hunt theorem  for functions on the polytorus $\mathbb{T}^\infty$. Obviously,  given some $f \in H_p(\mathbb{T}^\infty)$ with $1 < p < \infty$,
    there in general
    is no $z \in \mathbb{T}^\infty$ such that   $f(z) = \sum_{\alpha} \hat{f}(\alpha) z^\alpha$.

\begin{Coro}
Let  $f \in H_p(\mathbb{T}^\infty)$ and $1 < p < \infty$.  Then  for almost all $z \in \mathbb{T}^\infty$
\[
f(z) = \lim_N \sum_{\mathfrak{p}^\alpha \leq N} \hat{f}(\alpha) z^\alpha\,,
\]
and moreover,
\[
\bigg\| \sup_N\Big| \sum_{\mathfrak{p}^\alpha \leq N} \hat{f}(\alpha) z^\alpha \Big|^p \bigg\|_p
\leq C(p) \|f\|_p\,.
\]
\end{Coro}

\section{Vector-valued aspects of general Dirichlet series}
We now briefly discuss some results on vector-valued general Dirichlet series $\sum a_{n}e^{-\lambda_{n}s}$, where $(a_{n}) \subset X$ for some Banach space $X$. All results are related to the above topics on general  Dirichlet series
with scalar coefficients, and most of them are taken from \cite{vectorvalued}.

\medskip

In the ordinary case $\lambda = (\log n)$ the study of functional analytic aspects of vector-valued Dirichlet series
accumulated quite some amount of research, see e.g. \cite{CarandoDefantSevilla}, \cite{CarandoDefantSevilla1},
\cite{CarandoMaceraScottiTradacete},
\cite{CastilloMedinaGarciaMaestre},
\cite{DefantGarciaMaestrePerez},
\cite{DefantPopaSchwarting}, \cite{DefantSchwartingSevilla},
or \cite{DefantSevilla}. As an example we recall the following strong extension of \eqref{S12} from \cite{DefantGarciaMaestrePerez}: The maximal width of the strip of uniform, non absolute convergence of $X$-valued ordinary Dirichlet series
is given by the formula
\[
S(X) = 1 - \frac{1}{\text{cot}(X)}\,,
\]
where $\text{cot}(X)$ denotes the infimum over all $2 \leq q \leq \infty$ for which $X$
has cotype $q$ (see also \cite[Theorem 26.4]{Defant}).
\medskip

Given a Banach space $X$ and a frequency $\lambda$, we  denote by $\mathcal{D}(\lambda, X)$ the space of all $\lambda$-Dirichlet series with coefficients in $X$, and by   $\mathcal{D}_{\infty}(\lambda,X)$  the space of all $D\in \mathcal{D}(\lambda,X)$ which on $[Re>0]$ define bounded (and then necessarily holomorphic) functions. Moreover,
the by now obvious  definition of the a priori larger space $\mathcal{D}_{\infty}^{ext}(\lambda,X)$
is then clear.

\medskip

The Hahn-Banach theorem transports many results from scalar-valued Dirichlet series to vector-valued ones. For instance this way we may easily deduce
from Corollary \ref{almost periodic}
that $\left(\mathcal{D}_{\infty}(\lambda,X), \|\cdot\|_{\infty}\right)$ is a normed space with
$$\sup_{n \in \N} \|a_{n}\|_{X} \le \|D\|_{\infty}$$
for all $D \in \mathcal{D}_{\infty}(\lambda,X)$. In \cite{vectorvalued} it is shown that   $\mathcal{D}_{\infty}(\lambda)$ is complete if and only if $\mathcal{D}_{\infty}(\lambda,X)$ is.

\medskip

Another straight forward application of the Hahn-Banach theorem yields that the vector-valued version of Theorem  \ref{keylemma} holds, and consequently
also
the qualitative versions \eqref{coff} and \eqref{ruecken} of Bohr's theorem  (where in all inequalities the absolute value of course has to be replaced by the norm in $X$).

\medskip

Given  $1\le p\le \infty$, $X$ a Banach space and  $\lambda$  a frequency,
the definition of Hardy spaces of $X$-valued Dirichlet series extends naturally to the vector-valued situation.
 As in the scalar case we for any $\lambda$-Dirichlet group $(G,\beta)$ define the Banach space
$$H^{\lambda}_{p}(G,X):=\left\{ f \in L_{p}(G,X) \mid
\text{$f: \hat{G} \to X$ is supported on all $h_{\log n}, n \in \mathbb{N}$} \right\},$$
and consider  the vector-valued  Bohr map
$$\Bcal \colon H_{p}^{\lambda}(G,X) \hookrightarrow \mathcal{D}(\lambda,X), ~~ f \mapsto D=\sum \hat{f}(h_{\lambda_{n}}) e^{-\lambda_{n}s}\,.$$
Then it again turns out that  the Hardy space
$$\Hcal_{p}(\lambda,X):=\Bcal(H_{p}^{\lambda}(G, X))$$ together with the norm $\|D\|_{p}:=\|f\|_{p}$
 of all $X$-valued $\lambda$-Dirichlet series, forms a Banach space which is independent of the chosen $\lambda$-Dirichlet group.

 \medskip

In this final section of our survey we focus on two  different topics from \cite{vectorvalued}.  In the first one we ask up to which  extend Theorem \ref{Dstar} transfers to the  vector-valued setting? More precisely, for which Banach spaces $X$ and frequencies $\lambda$ does the equality $\mathcal{D}_{\infty}(\lambda,X)=\mathcal{H}_{\infty}(\lambda,X)$ hold (with an coefficient preserving isometry). Secondly, we revisit the maximal inequalities from the Sections \ref{Helson} and \ref{Carlesonchapter} in the vector-valued case.
\medskip
\subsection{Coincidence of vector-valued $\pmb{H_\infty}$'s}
We recall that a Banach space $X$ has the analytic Radon Nikodym property (short ARNP) if every
$f \in H_{\infty}(\D,X)$ has radial limits in almost all  $w \in\mathbb{T}$.
This happens if and only if $H_{\infty}(\T,X)=H_{\infty}(\D,X)$ (via an isometry identifying Fourier and
monomial coefficients).
Then it is straight forward to see that $\mathcal{D}_{\infty}((n),X)= \mathcal{H}_\infty((n),X)$
if and only if $X$ has ARNP.

\medskip

The main result of \cite{AntonioDefant} (see also \cite[Theorem 24.17]{Defant}) shows that this result holds true  for ordinary Dirichlet series, so for the frequency $\lambda = (\log n)$, and in the following theorem we state  an extension of all this for general Dirichlet series which is proven in \cite{vectorvalued}.

\begin{Theo} \label{ARNP}
Let $\lambda$ be any frequency and $X$ a non trivial Banach space with ARNP. Then
$\mathcal{D}_{\infty}(\lambda)=\Hcal_{\infty}(\lambda)$ if and only if
\begin{equation} \label{lufthansa}
\mathcal{D}_{\infty}(\lambda,X)=\Hcal_{\infty}(\lambda,X)\,.
\end{equation}
  In particular,
\eqref{lufthansa} holds whenever $X$ has ARNP and $\lambda$ satisfies $(BC)$.
\end{Theo}
But in fact there are Banach spaces $X$ without ARNP as well as  frequencies
$\lambda$
for which  $\mathcal{D}_{\infty}(\lambda)=\Hcal_{\infty}(\lambda)$, such that \eqref{lufthansa} is no equality. To see this we extend Theorem \ref{independent} to the vector-valued case. Denote by $\ell_{1}^{w}(X)$ the Banach space of all weakly summable $X$-valued sequences $(a_{n})$ with norm
$$w((a_{n}))=\sup_{x^{\prime} \in B_{X^{\prime}} } \sum_{n=1}^{\infty} |x^{\prime}(a_{n})|\,,$$
and by $\ell^{w,0}_{1}(X)$ the closed subspace of $\ell_{1}^{w}(X)$ consisting of all sequences $(a_{n})$
which are unconditionally summable.
Recall that  $\ell^{w,0}_{1}(X)=\ell^{w}_{1}(X)$ if and only if $X$ does not contain an isomorphic
copy of $c_{0}$. The following theorem is the announced extension of Theorem \ref{independent}  from \cite{vectorvalued}.
\begin{Theo} \label{independentvector} If $\lambda$ is $\Z$-linearly independent, then for every Banach space $X$ the following isometric inclusions hold:
$$\ell^{w,0}_{1}(X) \subset \Hcal_{\infty}(\lambda,X) \subset \mathcal{D}_{\infty}(\lambda,X) \subset\ell^{w}_{1}(X)$$
 In particular, $\mathcal{D}_{\infty}(\lambda,X)=\Hcal_{\infty}(\lambda,X)$ holds isometrically whenever  $c_{0}$ is
 no isomorphic subspace of $X$.
\end{Theo}
Banach lattices have ARNP if and only if  they don't contain an isomorphic copy of $c_0$. But in contrast to this, there are Banach spaces which fail to have  ARNP as well as an isomorphic copy of $c_0$ (see \cite[\S4]{Bukvalov}). Consequently, Theorem \ref{independentvector} shows that if $\lambda$ is $\Z$-linearly independent, then the equality $\mathcal{H}_{\infty}(\lambda,X)=\mathcal{D}_{\infty}(\lambda,X)$ does not necessarily imply that $X$ has ARNP.
\medskip

\subsection{Maximal inequalities - vector-valued}
Now we discuss the inequalities from the Sections \ref{Schauderbasischapter}, \ref{Helson} and \ref{Carlesonchapter}.
Recall that our proof of Theorem \ref{schauderbasisinHp} relies explicitely on the boundedness of the Riesz projection on $L_{p}(G)$, where $G$ is some Dirichlet group. This fact remains valid for $X$-valued functions
whenever $X$ is a so-called UMD-space which, as observed in \cite{vectorvalued}, is an immediate consequence of results from \cite{AsmarKellyMonty}.
\begin{Theo} Let $G$ be a Dirichlet group, $X$ a Banach space and $1<p<\infty$. Then $X$ has UMD if and only if the Riesz projection
$$R \colon L_{p}(G,X) \to L_{p}(G,X), ~~ f \mapsto \sum_{x\ge 0} \widehat{f}(h_{x}) h_{x}$$  is bounded.
\end{Theo}

As done in \cite{DefantSchoolmann} for the scalar case, a standard argument  leads to the following consequence
generalizing Theorem \ref{schauderbasisinHp} (see once again  \cite{vectorvalued}).

\begin{Coro}
For every frequency $\lambda$ and every UMD-space $X$ we have
$$\sup_N \Big\|S_N^{p,X}: \mathcal{H}_p(\lambda,X) \to \mathcal{H}_p(\lambda,X),\,\,
\sum a_n e^{-\lambda_n s} \mapsto \sum_{n=1}^N a_n e^{-\lambda_n s}\Big\| < \infty.
$$
%
%
%
\end{Coro}

Let us finally again turn to the  Carleson-Hunt theorem.
The proof of Theorem \ref{Carleson} given in \cite{DefantSchoolmann2}, extends in the following sense  from scalar to vector-valued functions.

\begin{Theo} \label{werderbr}
Assume that $X$ is a Banach space for which the Carleson-Hunt maximal inequality holds, i.e. for every $1 < p < \infty$ there
is some $C >0$ such that for every $f \in L_p(\mathbb{T},X)$
\begin{equation} \label{scalar}
  \bigg\| \sup_N \Big\| \sum_{n=-N}^N \hat{f}(n)  z^n  \Big\| \bigg\|_{L_p(\mathbb{T})}
  \leq C \|f\|_{L_p(\mathbb{T},X)}\,.
 \end{equation}
 Then for every frequency  $\lambda$ of integer type,
 every $\lambda$-Dirichlet group $(G, \beta)$,  and every
  $1 < p < \infty$ we have that for every $f \in H_p^\lambda(G,X)$
  $$f = \sum_{n=1}^\infty  \hat{f}(h_{\lambda_n})  h_{\lambda_n}  $$
  almost everywhere on $G$, and
   \begin{equation} \label{vector}
  \bigg\| \sup_N \Big\| \sum_{n=1}^N \hat{f}(h_{\lambda_n})  h_{\lambda_n}  \Big\| \bigg\|_{p}
  \leq C \|f\|_{p}\,,
 \end{equation}
where $C = C(p)$ is a constant which only depends on $p$.
\end{Theo}

In particular, if $X$ is a Banach lattice with  UMD, then  \eqref{scalar} holds true (see e.g. \cite{Francia}),
hence  in this case its $X$-valued variant \eqref{vector} is valid.
Actually, it is proven in \cite{vectorvalued} that for $\mathbb{Z}$-linearly independent frequencies $\lambda$  the assumption  UMD on $X$  in fact is  superfluous.
\begin{Prop}If $\lambda$ is $\mathbb{Z}$-linearly independent, then \eqref{vector} holds for all Banach spaces $X$.
\end{Prop}

Since the proof of Theorem \ref{maximalineq} needs for example tools like the
Hausdorff-Young inequality, there so far seems no  satisfying answer to the following final problem.

\begin{Prob} \label{Problem911} Under which assumptions on the Banach space $X$ and the frequency $\lambda$  for all $f\in H_{p}^{\lambda}(G,X)$ and $u>0$
\begin{equation*} \label{fr}
\left(\int_{G} \sup_{N} \left\|\sum_{n=1}^{N}\widehat{f}(h_{\lambda_{n}}) e^{-u\lambda_{n}}h_{\lambda_{n}}(\omega) \right\|_{X}^{p} d\omega \right)^{\frac{1}{p}} \le C \|f\|_{p}~?
\end{equation*}
\end{Prob}
If there is such an estimate, one may ask on which of the parameters $\lambda, X, p$ and $u$ the constant $C$ actually depends.
Recall that in Theorem \ref{maximalineq}, which deals with the scalar case $X=\mathbb{C}$ and $\lambda$'s with $(BC)$, the constant $C$ in $(\ref{ineqmaxi})$ does not depend on $p$ (but on $u$ and $\lambda$).
 This is a useful fact since by  Remark \ref{ptoinfity}, Theorem \ref{maximalineq} then implies  Bohr's theorem for  $\lambda$'s with $(BC)$. Since Bohr's theorem in the vector-valued case holds for every Banach space $X$ provided that $\lambda$ has $(LC)$, Problem \ref{Problem911} may hold for every Banach space $X$ and $\lambda$'s satisfying $(LC)$ with some constant $C$   independent of $p$.
 If we ask for a constant $C$ depending on $p$, then in  view of Theorem \ref{werderbr} it may happen that Problem \ref{Problem911} has an affirmative answer 
 for every frequency $\lambda$ but with restrictions on the Banach space $X$.

\end{document}